\documentclass[11pt,a4paper]{article}%
\usepackage{graphicx}
\usepackage{amsmath}
\usepackage{amsfonts}
\usepackage{amssymb}%
\newtheorem{theorem}{Theorem}

\newtheorem{corollary}[theorem]{Corollary}

\newtheorem{definition}[theorem]{Definition}

\newtheorem{lemma}[theorem]{Lemma}

\newtheorem{proposition}[theorem]{Proposition}

\newenvironment{proof}[1][Proof]{\textbf{#1.} }{\ \rule{0.5em}{0.5em}}

\begin{document}

\title{Perturbed Markov Chains}
\author{Eilon Solan\thanks{MEDS Department, Kellogg School of Management, Northwestern
University, \emph{and} School of Mathematical Sciences, Tel Aviv University,
Tel Aviv 69978, Israel. e-mail: eilons@post.tau.ac.il} \ and Nicolas
Vieille\thanks{D\'{e}partement Finance et Economie, HEC, 1, rue de la
Lib\'{e}ration, 78 351 Jouy-en-Josas, France. e-mail: vieille@hec.fr}%
\ \thanks{We acknowledge the financial support of the Arc-en-Ciel/Keshet
program for 2001/2002. The research of the second author was supported by the
Israel Science Foundation (grant No. 03620191).}}
\maketitle

\begin{abstract}
We study irreducible time-homogenous Markov chains with finite state space in
discrete time. We obtain results on the sensitivity of the stationary
distribution and other statistical quantities with respect to perturbations of
the transition matrix. We define a new closeness relation between transition
matrices, and use graph-theoretic techniques, in contrast with the matrix
analysis techniques previously used.

\end{abstract}

\bigskip


\bigskip

Primary subject classification: 60J10

Secondary subject classification: 60F10

\newpage

\section{Introduction}

The present paper concerns irreducible time-homogenous Markov chains with a
finite state space in discrete time. We are interested in the effects of
perturbations of the transition matrix on the stationary distribution and on
other statistical quantities.

This question has a long history, starting with Schweizer \cite{Sch68}. Let
$S$ be a finite set, and $q=(q(t\mid s))_{s,t\in S}$ an irreducible transition
matrix. Let $\widehat{q}=((\widehat{q}(t\mid s))_{s,t\in S}$ be another
transition matrix. Denote by $\mu=(\mu_{s})_{s\in S}$ and $\widehat{\mu
}=(\widehat{\mu}_{s})_{s\in S}$ the stationary distributions that correspond
to $q$ and $\widehat{q}$ respectively. Schweizer \cite{Sch68} estimated
$\Vert\mu-\widehat{\mu}\Vert_{\infty}$ as a function of $\Vert q-\widehat
{q}\Vert_{\infty}$, using the fundamental matrix of Kemeny and Snell
\cite{KemSne60}. A vast literature has followed up Schweizer, and provided
various estimates to $\Vert\mu-\widehat{\mu}\Vert_{\infty}$ and $\left\vert
\frac{\mu_{s}-\widehat{\mu}_{s}}{\mu_{s}}\right\vert $ ($s\in S$) as a
function of $\Vert q-\widehat{q}\Vert_{\infty}$. See, e.g., \cite{Mey75},
\cite{Sen88}, \cite{Sen93}, \cite{Les94}, \cite{KNS98}, \cite{ChoMey00}.

O'Cinneide \cite{Oci93} studied the effects of entry-wise relative
perturbations on the stationary distribution; that is, he provided a bound on
$\left|  \frac{\mu_{s} - \widehat\mu_{s}}{\mu_{s}} \right| $ as a function of
$\max_{s,t \in S} \left\{  \frac{q(t \mid s)}{\widehat q(t \mid s)},
\frac{\widehat q(t \mid s)}{ q(t \mid s)}\right\} $.

Our paper differs from the existing literature in three respects. First, we
introduce a new way to measure the difference between two transition matrices.
In the spirit of the entry-wise relative perturbations of O'Cinneide
\cite{Oci93}, our measure is $\max\left\{  \left|  1 - \frac{\widehat q(t \mid
s)}{ q(t \mid s)}\right|  \right\} $, but the maximum is not taken over
\emph{all} pairs $s,t \in S$. Rather, the maximum is taken over all
\emph{frequent} transitions $s \to t$. Formally, the notion of closeness
between $q$ and $\widehat{q}$ is the following. For a transition matrix $q$,
we define $\zeta_{q} = \min_{\emptyset\subset C\subset S}\sum_{s\in C}\mu
_{s}q(\overline{C}\mid s)$, where $\overline{C} = S \setminus C$ is the
complement of $C$ in $S$. This is a variant of the conductance (see, e.g.,
\cite{JerSin89}, \cite{LasRav01}, \cite{LovSim90}). Given $\varepsilon
,\beta>0$, we say that $\widehat{q}$ is $(\varepsilon,\beta)$-close to $q$ if
for every two states $s,t\in S$, $\left|  1-\frac{\widehat{q}(t|s)}%
{q(t|s)}\right|  \leq\beta$ whenever (a) $\mu_{s}q(t|s)\geq\varepsilon
\zeta_{q}$ or (b) $\mu_{s}\widehat{q}(t|s)\geq\varepsilon\zeta_{q}$. Condition
(a) holds whenever the transition from $s$ to $t$ occurs frequently. Condition
(b) is not analogous to (a), since it involves the stationary distribution of
$q$, and the transition matrix $\widehat{q}$.

Provided $\varepsilon$ and $\beta$ are small enough, we show that
$\widehat {q}$ is irreducible, and
$\frac{\mu_{s}}{\widehat{\mu}_{s}}$ is close to 1 for each $s \in S$.

The logic behind this closeness relation, and its novelty, is that even large
perturbations that occur in rarely visited states should not affect too much
the stationary distribution. This point is illustrated in the example that is
studied in the next section.

A motivation for using this closeness measure is the following statistical
implication. Assume a statistician observes a realization $s_{1},...,s_{N}$
of (the first $N$ components of) a Markov chain with unknown irreducible
transition matrix $q$. Given these observations, he computes the empirical
transition matrix $\widehat{q}$ and the invariant distribution $\widehat{\mu}$
of $\widehat{q}$ (if $s_{N}=s_{1}$, $\widehat{\mu}$ coincides with the
empirical occupancy measure). For fixed $\varepsilon,\beta > 0$, if $N$ is large enough, $\widehat{q}$ is
$(\varepsilon,\beta)$-close to $q$ with probability close to one, hence
$\widehat{\mu}$ is an accurate estimate of the invariant distribution of the
underlying Markov chain.

Second, instead of using matrix analysis, we use the graph techniques
developed by Freidlin and Wenzell \cite{FreWen} and extensively used in the
analysis of Markov chains with rare transitions (see, e.g., Catoni
\cite{Cat99}).

Third, our approach allows us to estimate the sensitivity of other statistical
quantities, such as the exit distribution from a given set and the average
length of visits to a given set.

The paper is organized as follows. Section \ref{secmodele} contains the
statements of the main results. We also study there an example that
illustrates the advantage of our closeness relation, and compares the bound we
derive to existing bounds. Section \ref{sec graphes} briefly recalls standard
formulas, and states few elementary properties of graphs. Section
\ref{sec proofmain} is devoted to the proof of the main result. The last
section deals with a variation of the main result.

\section{Notations and results\label{secmodele}}

Let $S$ be a finite set, fixed through the paper, with at least two elements.
For every subset $C\subseteq S$, $\overline{C}=S\setminus C$ is the complement
of $C$ in $S$, $|C|$ is its cardinality, and $\Delta(C)$ is the set of
probability distributions over $C$. For $s \in C \subseteq S$, we write $C
\setminus s$ instead of the more cumbersome $C \setminus\{s\}$.

\subsection{Main result}

Let $q$ be an irreducible transition matrix over $S$, with stationary
distribution $\mu=(\mu_{s})_{s\in S}$. For every $C \subseteq S$ we denote
$\mu_{C} = \sum_{s \in C} \mu_{s}$. Let $\widehat{q}$ be another transition
matrix over $S$. Assuming $\widehat{q}$ is irreducible, we wish to bound the
distance between $\mu$ and $\widehat{\mu}$.

Our notion of closeness of $\widehat{q}$ to $q$ involves a measure of how
mixing $q$ is$.$ Our measure involves the quantity%
\begin{equation}
\zeta_{q}=\min_{\emptyset\subset C\subset S}\sum_{s\in C}\mu_{s}q(\overline
{C}\mid s), \label{zeta}%
\end{equation}
which is a variant of the \emph{conductance}, see e.g. \cite{JerSin89},
\cite{LasRav01}, \cite{LovSim90}. Given $C\subset S$, the quantity $\sum_{s\in
C}\mu_{s}q(\overline{C}\mid s)\,\;$measures the average frequency of
transitions out of $C$. Hence, $\zeta_{q}\,$, being the lowest such frequency,
is a measure of how isolated a subset $C$ may be. Formally,

\begin{definition}
Let $\varepsilon,\beta>0$. We say that a transition matrix $\widehat{q}$ is
\emph{$(\varepsilon,\beta)$-close to $q$} if for every two states $s,t\in S$,
$\left|  1-\frac{\widehat{q}(t|s)}{q(t|s)}\right|  \leq\beta$ whenever (a)
$\mu_{s}q(t|s)\geq\varepsilon\zeta_{q}$ or (b) $\mu_{s}\widehat{q}%
(t|s)\geq\varepsilon\zeta_{q}$.
\end{definition}

Note that this closeness notion is not symmetric, since we use only the
stationary distribution of $q$.

Denote $L=\sum_{n=1}^{|S|-1}\binom{|S|}{n}n^{|S|}$. We now state our main result.

\begin{theorem}
\label{theoremmain} Let $\beta\in(0,1/2^{|S|})$ and let $\varepsilon
\in(0,\frac{\beta(1-\beta)}{L\times|S|^{4}})$. For every irreducible
transition matrix $q$ on $S$ and every transition matrix $\widehat{q}$ that is
$(\varepsilon,\beta)$-close to $q$:

\begin{enumerate}
\item $\widehat q$ is irreducible.

\item Its stationary distribution $\widehat{\mu}$ satisfies $\left|
1-\frac{\widehat{\mu}_{s}}{\mu_{s}}\right|  \leq18\beta L$ for each $s\in S$.
\end{enumerate}
\end{theorem}

\subsection{An example and comparison with existing bounds}

As mentioned in the introduction, many authors provided bounds for the
sensitivity of the stationary distribution.

We now study an example that first highlights the logic behind the closeness
relation, and second shows that in some cases, the bound we give is better
than existing bounds.

Fix $\delta\in(0,1/2)$. Take a Markov chain with three states and transition
matrix as follows.
\[
q(2\mid1)=1-\delta,\ q(3\mid1)=\delta,\ q(1\mid2)=q(1\mid3)=1.
\]

Thus, in every other stage the process visits state 1. In particular, the
stationary distribution $\mu$ is given by
\[
\mu_{1} = \frac{1}{2}, \ \mu_{2} = \frac{1-\delta}{2}, \ \mu_{3} =
\frac{\delta}{2}.
\]
The quantity $\zeta_{q}$ is given by
\[
\zeta_{q} = \min\left\{ \frac{1}{2}, \frac{1-\delta}{2}, \frac{\delta}{2},
\frac{1}{2} \times\delta, \frac{1}{2} \times(1-\delta), \frac{1-\delta}{2} +
\frac{\delta}{2} \right\}  = \frac{\delta}{2}.
\]

Fix $0 < \beta< 1/8 = 1/2^{3} $, $0 < \varepsilon< \beta(1-\beta)/3^{4}L$ and
$0 < \eta< \varepsilon$. In particular, one can choose $\eta= \beta/C$, where
$C > 1$ is some fixed scalar, independent of $\delta$.

Define a transition matrix $\widehat q$ by
\[
\widehat q(2 \mid1) = 1-\delta, \ \widehat q(3 \mid1) = \delta, \ \widehat q(1
\mid2) = 1, \ \widehat q(1 \mid3) = 1-\eta, \ \widehat q(2 \mid3) = \eta.
\]
Thus, we only change the transitions out of state 3. Moreover, whereas the
relative size of the transition $3 \to1$ changed moderately, the relative size
of the transition $3 \to2$ changed dramatically.

We first argue that $\widehat q$ is $(\varepsilon,\beta)$-close to $q$.
Indeed, $\left|  1 - \frac{\widehat q(t \mid s)}{q(t \mid s)} \right|  = 0$
whenever $s \neq3$, and $\left|  1 - \frac{\widehat q(1 \mid3)}{q(1 \mid3)}
\right|  = \eta< \beta$. The claim now follows, since $\mu_{3} q(2 \mid3) = 0$
and $\mu_{3} \widehat q(2 \mid3) = \frac{\delta}{2} \times\eta< \zeta_{q}
\times\varepsilon$.

Thus, Theorem \ref{theoremmain} states that $\left|  \frac{\mu_{3} -
\widehat\mu_{3}}{\mu_{3}} \right|  \leq18L \beta= 18LC \eta$.

This example highlights the logic behind our closeness relation. Since state 3
is rarely visited, the stationary distribution is not too sensitive to changes
in transitions out of this state, even though these changes are relatively large.

\bigskip

Cho and Meyer \cite{ChoMey00} bound the sensitivity of the stationary
distribution by the mean first passage time:
\begin{equation}
\left\vert \frac{\mu_{s}-\widehat{\mu}_{s}}{\mu_{s}}\right\vert \leq
\frac{\Vert q-\widehat{q}\Vert_{\infty}}{2}\times\max_{t\neq s}M_{t,s}%
,\label{chomey}%
\end{equation}
where $M_{t,s}$ is the mean first passage time from state $t$ to state $s$.

Observe that in the example $M_{2,3}=2/\delta$ while $M_{1,3}=2/\delta-1$,
hence the bound provided in \cite{ChoMey00} for $\left\vert \frac{\mu
_{3}-\widehat{\mu}_{3}}{\mu_{3}}\right\vert $ is $\eta\times\frac{1}{\delta}$,
which is worse than our bound when $\delta$ is small.

In the case of two-state chains, the bound (\ref{chomey}) is very close to the
bound that can be derived from Theorem \ref{theoremmain}, up to a universal
constant. Indeed, (\ref{chomey}) then reduces to
\begin{equation}
\left\vert \frac{\mu_{s}-\widehat{\mu}_{s}}{\mu_{s}}\right\vert \leq\frac
{1}{2}\max\left\{\left\vert q(2|1)-\widehat{q}(2|1)\right\vert ,\left\vert
q(1|2)-\widehat{q}(1|2)\right\vert \right\}\times\max\left\{\frac{1}{q(2|1)},\frac
{1}{q(1|2)}\right\}.\label{chomey2}%
\end{equation}
On the other hand, $\zeta_{q}=\min\left\{  \mu_{1}q(2|1),\mu_{2}q(1|2)\right\}$.
Hence, Theorem \ref{theoremmain} yields%
\[
\left\vert \frac{\mu_{s}-\widehat{\mu}_{s}}{\mu_{s}}\right\vert \leq
36\beta\text{, with }\beta=\max\left\{  \left\vert \frac{q(2|1)-\widehat
{q}(2|1)}{q(2|1)}\right\vert ,\left\vert \frac{q(1|2)-\widehat{q}%
(1|2)}{q(1|2)}\right\vert \right\}  ,
\]
provided $\beta < 1/4$ - an inequality which is slightly more precise than (\ref{chomey2}), up to a
constant 72.

\bigskip

Kirkland et al \cite{KNS98} bound the relative sensitivity by
\[
\left\vert \frac{\mu_{s}-\widehat{\mu}_{s}}{\mu_{s}}\right\vert \leq\frac
{1}{2}\frac{\Vert q-\widehat{q}\Vert_{\infty}}{2}\times\Vert A_{s}^{-1}%
\Vert_{\infty},
\]
where $A$ is the fundamental matrix of \cite{KemSne60}, and $A_{s}$ is the
$(n-1)\times(n-1)$ submatrix of $A$ obtained by deleting the $s$'th row and column.

One can verify that in our example $A_{3} = \left(
\begin{array}
[c]{cc}%
1 & -1+\delta\\
-1 & 1
\end{array}
\right) $, so that $A^{-1}_{3} = \left(
\begin{array}
[c]{cc}%
1/\delta & -1+1/\delta\\
1/\delta & 1/\delta
\end{array}
\right) $. Thus, the bound for $\left|  \frac{\mu_{3} - \widehat\mu_{3}}%
{\mu_{3}} \right| $ is $\eta/2\delta$, which is worse than our bound when
$\delta$ is small.

\bigskip

Note that the entry-wise ratio bound given by O'Cinneide \cite{Oci93} is not
useful in this example, since $q(2 \mid3) = 0$ while $\widehat q(2 \mid3) > 0$.

\subsection{The transition matrix changes in a subset of states}

In some cases, the transition matrix is perturbed only in a subset $S_{1}
\subset S$, and restricted to $S_{1}$ the transition matrix is sufficiently
mixing, in the sense that the probability to reach any state in $S_{1}$ before
leaving $S_{1}$ is bounded from below. In this case, instead of taking the
conductance in the definition of the closeness relation, one can take another
quantity, which is, in a sense, the conductance restricted to $S_{1}$.

Such a case occurs, for example, if the state space can be partitioned into
some subsets, the transition matrix is mixing in each subset while the
probability to move from one subset to another is small, and one perturbs the
transition matrix only in one of the subsets.

Let $S_{1}$ be a subset of $S,$ with $\left|  S_{1}\right|  >1$. Define
\[
\zeta_{q}^{1}=\min_{\emptyset\neq C\subset S_{1}}\sum_{s\in C}\mu
_{s}q(\overline{C}\mid s).
\]
Let $(\mathbf{s}_{n})$ be a Markov chain with transition matrix $q$. We denote
by $\mathbf{P}_{s,q}$ the law of $(\mathbf{s}_{n})$ when the initial state is
$s$, and by $\mathbf{E}_{s,q}$ the corresponding expectation.

For every proper subset $C$ of $S$ we let $T_{C}=\min\left\{  n\geq0,
\mathbf{s}_{n}\in C\right\}  $ denote the first hitting time of $C$ and
$T_{C}^{+}=\min\left\{  n\geq1, \mathbf{s}_{n}\in C\right\}  $ the first
return to $C$ . By convention, the minimum over an empty set is $+\infty$.

\begin{definition}
Let $\varepsilon,\beta>0$. We say that a transition matrix $\widehat{q}$ is
\emph{$(\varepsilon,\beta)$-close to $q$ on $S_{1}$} if for every two states
$s,t\in S$, $\left|  1-\frac{\widehat{q}(t|s)}{q(t|s)}\right|  \leq\beta$
whenever (a) $\mu_{s}q(t|s)\geq\varepsilon\zeta_{q}^{1}$ or (b) $\mu
_{s}\widehat{q}(t|s)\geq\varepsilon\zeta_{q}^{1}$.
\end{definition}

We now state the Theorem that corresponds to Theorem \ref{theoremmain}.

\begin{theorem}
\label{theoremvariation}Let $\beta\in(0,1/2^{|S|})$, $a>0$ and $\varepsilon
\in(0,\frac{1}{2}\left(  \frac{a}{L}\right)  ^{\left|  S\right|  }\times
\frac{\beta(1-\beta)}{L\times|S|^{4}})$. Let $q$ be an irreducible transition
matrix such that $\mathbf{P}_{s,q}(T_{\overline{S}_{1}\cup\{t\}}^{+}%
=T_{\{t\}}^{+})\geq a$ for every $s,t\in S_{1}$. Then, for every transition
matrix $\widehat{q}$ that is $(\varepsilon,\beta)$-close to $q$ on $S_{1}$ and
that coincides with $q$ on $S\backslash S_{1}$, we have

\begin{enumerate}
\item \label{variantq} All states of $S_{1}$ belong to the same recurrent set
$R$ for $\widehat{q}$.

\item \label{variantmu} The stationary distribution $\widehat{\mu}$ of
$\widehat{q}$ on $R$ satisfies
\begin{equation}
\left|  1-\frac{\widehat{\mu}(s|S_{1})}{\mu(s|S_{1})}\right|  \leq18\beta
L\text{, for each }s\in S_{1},
\end{equation}
where $\mu(s \mid S_{1}) = \mu_{s} / \mu_{S_{1}}$.
\end{enumerate}
\end{theorem}

\bigskip

Note that the claims in Theorem \ref{theoremvariation} differ from those in
Theorem \ref{theoremmain}. It is no longer claimed that $\widehat{q}$ is
irreducible, nor that the unconditional stationary distributions $\mu$ and
$\widehat{\mu}$ are close. The statements in Theorem \ref{theoremvariation}
are optimal in this respect. This is due to the fact that the quantity
$\zeta_{q}^{1}$ contains no information on the frequency of transitions out of
$S_{1}$. To emphasize this point, consider the following example.

Assume that $S=\left\{  a,b,c\right\}  $ and $S_{1}=\left\{  a,b\right\}  $.
Let $\varepsilon,\beta\in(0,1/2)$ be given. Let two additional parameters
$\lambda$ and $\eta$ be given in $(0,1)$, and define $q$ as follows. From
state $a$ (resp. $b$) a chain with transition matrix $q$ moves to $c$ with
probability $\eta$, and otherwise to $b$ (resp. to $a$). From state $c,$ the
chain remains in $c$ with probability $1-\lambda$, and otherwise moves to $a$
or $b$ with equal probability $\frac{1}{2}\lambda$.

Plainly, $q$ is irreducible, and the value of $\mu_{a}=\mu_{b}$ depends on the
ratio $\lambda/\eta$: this common value may be arbitrary close to $0$ (resp.
to $1/2$) provided $\lambda/\eta$ is close enough to $0$ (resp. to $+\infty$).
Note that $\zeta_{q}^{1}=\mu_{a}q(\left\{  b,c\right\}  |a)=\mu_{a}$. Let now
$\widehat{q}$ be defined exactly as $q$, except that the parameter $\eta$ is
replaced by another parameter $\widehat{\eta}\in\left[  0,1\right]  $. As soon
as $\eta,\widehat{\eta}<\min\{\varepsilon,\beta\}$, $\widehat{q}$ is
($\varepsilon,\beta)$-close to $q$. This is in particular the case if
$\widehat{\eta}=0$, in which case $\widehat{q}$ fails to be irreducible. On
the other hand, even if $\widehat{\eta}>0$, the values of $\eta,\widehat{\eta
}$ and $\lambda$ can be chosen in such a way that the inequalities
$\eta,\widehat{\eta}<\min\{\varepsilon,\beta\}$ are satisfied, and $\eta
\ll\lambda\ll\widehat{\eta}$.\ Hence, even if $\widehat{q}$ is irreducible,
its unconditional stationary distribution $\widehat{\mu}$ may be arbitrarily
far from $\mu$.

\bigskip

\subsection{Sensitivity of other quantities}

Our graph-theoretic approach allows us to obtain information on other
quantities of interest. We here present the statements of the corresponding results.

We let $\mathbf{Q}_{s,q}(\cdot|C)$ denote the law of the exit state from $C$:
$\mathbf{Q}_{s,q}(t|C)=\mathbf{P}_{s,q}(T_{\overline{C}}=T_{t})$ for $t\notin
C$. Next, we set
\begin{align}
\nu_{C}(s)  &  :=\frac{\sum_{t\in\overline{C}}\mu_{t}q(s|t)}{\sum
_{t\in\overline{C}}\mu_{t}q(C|t)}\text{ for }C\subset S\text{ and }s\in C,
\hbox{ and}\label{averageentry}\\
K_{C}  &  :=\sum_{s\in C}\nu_{C}(s){\mathbf{E}}_{s,q}[e_{C}]\text{ for
}C\subset S.\nonumber
\end{align}
The numerator (resp. the denominator) in (\ref{averageentry}) is the long run
frequency of transitions from $\overline{C}$ to $s$ (resp. from $\overline{C}$
to $C$). Thus, $\nu_{C}^{{}}(s)$ is the probability that the first stage in
$C$ the process visits is $s$, while $K_{C}$ is the average length of a visit
to $C$.

Assuming $\widehat{q}$ is irreducible, the corresponding quantities for
$\widehat{q}$ will be denoted by ${\mathbf{Q}}_{s,\widehat{q}}$, $\widehat
{\nu}_{C}(s)$ and $\widehat{K}_{C}$. We now state the results on
$\mathbf{Q}_{s,\widehat{q}}$ and $\widehat{K}_{C}$ that hold in the framework
of Theorems \ref{theoremmain} and \ref{theoremvariation} respectively.

\bigskip

\begin{theorem}
\label{theorem main2}Set $c=2\left|  S\right|  ^{2}$. Under the assumptions of
Theorem \ref{theoremmain}, the following holds: for each $C\subset S$,

\begin{enumerate}
\item $\left\|  {\mathbf{Q}}_{s,q}(\cdot\mid C)-{\mathbf{Q}}_{s,\widehat{q}%
}(\cdot\mid C)\right\|  <12\beta L$ for every $s\in C$

\item $\frac{1}{c}K_{C}\leq\widehat{K}_{C}\leq cK_{C}$.
\end{enumerate}
\end{theorem}

\begin{theorem}
\label{theorem variation2}Set $c=2\left|  S\right|  ^{2}$. Under the
assumptions of Theorem \ref{theoremvariation}, the following holds: for each
$C$ $\subset S_{1}$,

\begin{enumerate}
\item $\left\|  {\mathbf{Q}}_{s,q}(\cdot\mid C)-{\mathbf{Q}}_{s,\widehat{q}%
}(\cdot\mid C)\right\|  <12\beta L$ for every $s\in C$

\item $\frac{1}{c}K_{C}\leq\widehat{K}_{C}\leq cK_{C}$

\item $\frac{1}{c}K_{S_{1}}\leq\widehat{K}_{S_{1}}\leq cK_{S_{1}}$ or
$K_{S_{1}},\widehat{K}_{S_{1}}\geq\frac{1}{2\varepsilon\left|  S\right|
}\times\frac{\mu_{S_{1\;}}}{\zeta_{q}^{1}}$.
\end{enumerate}
\end{theorem}

We let $q$ be an irreducible transition matrix over $S$. It is fixed
throughout the paper.

\section{Preliminaries\label{sec graphes}}


Our computations are based on formulas due to Freidlin and Wenzell
\cite{FreWen}, that express stationary distribution, exit distributions and
expected hitting times in graph-theoretic terms. For a discussion of some
applications, we refer to Catoni \cite{Cat99}. These tools have also been used
in the context of stochastic games in \cite{VieIJM2} and \cite{SolVie98}.

The weight of a graph is obtained from the transition probabilities
corresponding to the different edges of the graph. We recall these formulas in
section \ref{subsecreminder}. Next, we compare the weights of a given graph
under a transition matrix $\widehat{q}$ that is close to $q$.

\subsection{Reminder\label{subsecreminder}}

Given $C\subset S$, a $C$\emph{-graph} is a directed graph without cycle $g$
over $S$ such that:\footnote{Our $C$-graphs correspond to $\overline{C}%
$-graphs in \cite{FreWen}, \cite{Cat99}.}

\begin{itemize}
\item For $s\in C$, there is exactly one edge starting at $s,$ denoted by
$(s,g(s))$.

\item For $s\in\overline{C}$, there is no edge starting at $s$.
\end{itemize}

Thus, given $s\in C$, there is a unique path starting at $s$ and ending at
some $t\in\overline{C}$. We say that $s$ \emph{leads to }$t$ \emph{along} $g$.
We denote by $G(C)$ the set of $C$-graphs; for $s\in C,t\in\overline{C}$,
$G_{s,t}(C)$ is the subset of graphs $g\in G(C)$ such that $s$ leads to $t$
along $g$. Note that $G(C)$ depends only on $C$, and not on the transition
matrix. Note also that $L$ bounds the number of graphs: $L\geq\sum
_{\emptyset\subset C\subset S}|G(C)|$.

We identify each $C$-graph $g$ with the collection of its edges: $g=\cup_{s\in
C}\{(s,g(s))\}$.

Given $D\subseteq C$, and $g\in G(C)$, the \emph{restriction of }$g$ \emph{to
}$D$ is defined to be the subgraph of $g$ that contains exactly those edges of
$g$ that start in $D.$ Thus, it is the $D$-graph $g^{\prime}= \cup_{s\in
D}\left\{  (s,g(s))\right\}  $.

For every $g\in G(C)$, we define the \emph{weight} of $g$ under $q$ by
\[
p(g):=\prod_{(s,t)\in g}q(t|s)\text{.}%
\]

The following formulas were derived by Freidlin and Wentzell \cite{FreWen},
Lemmas 6.3.1, 6.3.4 and 6.3.3. For more direct statements and alternative
proofs see Catoni \cite{Cat99}.

\begin{proposition}
[Freidlin-Wenzell, 1984]\label{fw} Let $(S,q)$ be a Markov chain.

\begin{itemize}
\item If $q$ is irreducible then for every $s \in S$
\begin{equation}
\mu_{s} = \frac{\sum_{G(S\backslash\left\{  s\right\}  )}p(g)}{\sum_{y\in
S}\sum_{G(S\backslash\left\{  y\right\}  )}p(g)}. \label{invariant}%
\end{equation}

\item For every proper subset $C$ of $S$ and every $s\in C$,
\begin{equation}
\mathbf{E}_{s,q}\left[  T_{\overline{C}}\right]  =\frac{\sum_{G(C\backslash
\left\{  s\right\}  )}p(g)+\sum_{t\in C,t\neq s}\sum_{G_{s,t}(C\backslash
\left\{  t\right\}  )}p(g)}{\sum_{G(C)}p(g)}, \label{exit}%
\end{equation}
and
\begin{equation}
\mathbf{Q}_{s,q}(t|C)=\frac{\sum_{G(C)}p(g)}{\sum_{G_{s,t}(C)}p(g)}%
\hbox{ for each } t\notin C. \label{distribution}%
\end{equation}

\end{itemize}
\end{proposition}

\subsection{Basic properties\label{secbasic}}

In this section we provide basic properties of weights of graphs. The
transition matrix $q$ is here arbitrary.

\begin{definition}
Let $C$ be a proper subset of $S$, and let $\eta>0$. A graph $g\in G(C)$ is
\emph{$\eta$-maximal} if
\[
p(g)\geq\eta\max_{g^{\prime}\in G(C)}p(g^{\prime})\text{.}%
\]

\end{definition}

We denote by $G^{\eta}(C)$ the set of $\eta$-maximal $C$-graphs. For
simplicity of notations, we do not emphasize the dependency of $G^{\eta}(C)$
on the transition matrix. Clearly, $G^{\eta}(C)$ is non-empty, for every
$\eta\leq1$ and $C\subset S$. It is worth listing a few basic properties of
graphs that we use repeatedly.

\begin{proposition}
\label{prop graphs}

\begin{itemize}
\item[\textbf{P0}] Let $C_{1}\cap C_{2}=\emptyset$, and $g_{i}\in G(C_{i})$,
for $i=1,2$. If all paths of $g_{1}$ lead to $\overline{C_{1}\cup C_{2}}$,
then $g_{1}\cup g_{2}$ is a $C_{1}\cup C_{2}$-graph.

\item[\textbf{P1}] Let $C_{1}\cap C_{2}=\emptyset$, $g\in G^{\eta}(C_{1}\cup
C_{2})$, and $g_{i}$ the restriction of $g$ to $C_{i}$. If all paths of
$g_{2}$ lead to $\overline{C_{1}\cup C_{2}}$, then $g_{1}\in G^{\eta}(C_{1})$.

\item[\textbf{P2}] Let $C_{1}\cap C_{2}=\emptyset$, and $g_{i}\in G^{\eta_{i}%
}(C_{i})$ for $i=1,2$.\ If $g_{1}\cup g_{2}$ is a $C_{1}\cup C_{2}$-graph,
then it is $\eta_{1}\eta_{2}$-maximal.
\end{itemize}
\end{proposition}

\begin{proof}
\textbf{P0} and \textbf{P2} follow from the definitions. We now show that
\textbf{P1} holds. Otherwise, there is $g^{\prime}_{1} \in G(C_{1})$ such that
$p(g_{1}) < \eta p(g^{\prime}_{1})$. By \textbf{P0}, $g^{\prime}= g^{\prime
}_{1} \cup g_{2}$ is in $G(C_{1} \cup C_{2})$, but $p(g) < \eta p(g^{\prime}%
)$, a contradiction.
\end{proof}

\bigskip

Note that \textbf{P1} needs not hold without the condition that all paths of
$g_{2}$ lead to $\overline{C_{1}\cup C_{2}}$. Indeed, take $S = \{1,2,3,4\}$,
$C_{1} = \{1\}$, $C_{2} = \{2\}$, and $q(2\mid1) = q(1 \mid2) = 1-q(3 \mid1) =
1-q(4 \mid2) = 2/3$. The $C_{2}$-graph $g_{1} = (2 \to4)$ is $1/2$-maximal,
and the $C_{1} \cup C_{2}$-graph $(1 \to2, 2 \to4)$ is 1-maximal.

\begin{lemma}
\label{lemma h} Let $C$ be a proper subset of $S$, let $\eta>0$, and let $H$
be a set of graphs such that $G^{\eta}(C)\subseteq H\subseteq G(C)$. Then
\[
0\leq\frac{\sum_{g\in G(C)}p(g)}{\sum_{g\in H}p(g)}-1<\eta L.
\]
In particular,
\[
0\leq1-\frac{\sum_{g\in H}p(g)}{\sum_{g\in G(C)}p(g)}<\eta L.
\]

\end{lemma}

\begin{proof}
Since $H \subseteq G(C)$, and by the definition of $G^{\eta}(C)$,
\[
0 \leq\frac{\sum_{g \in G(C)} p(g)}{\sum_{g \in H} p(g)} - 1 = \frac{\sum_{g
\in G(C) \setminus H} p(g)}{\sum_{g \in H} p(g)} \leq\frac{\sum_{g \in G(C)
\setminus G^{\eta}(C)} p(g)}{\sum_{g \in G^{1}(C)} p(g)} < \eta L,
\]
as desired.
\end{proof}

\section{Proof of the main results\label{sec proofmain}}

We here prove Theorems \ref{theoremmain} and \ref{theorem main2}. We let
$\varepsilon,\beta\in(0,1)$ satisfy the assumptions of Theorem
\ref{theoremmain}, and $\widehat{q}$ be another transition matrix over $S$. We
assume that $\widehat{q}$ is $(\varepsilon,\beta)$-close to $q$.

\subsection{On graphs\label{secestimates}}

For every proper subset $C$ of $S$ and every $\eta>0$, we denote by
$\widehat{G}^{\eta}(C)$ the set of $\eta$-maximal graphs under $\widehat{q}$.
For every $C$-graph $g$, $\widehat{p}(g)=\prod_{s\in C}\widehat{q}(g(s)\mid
s)$ is the weight of $g$ under $\widehat{q}$.

\begin{lemma}
\label{lemmabasic1} For every proper subset $C$ of $S$,
\begin{equation}
\frac{1-\beta}{|S|^{2}}\sum_{s\in C}\mu_{s}{q}(\overline{C}|s)\leq\sum_{s\in
C}\mu_{s}\widehat{q}(\overline{C}|s)\leq(1+\beta)|S|^{2}\sum_{s\in C}\mu
_{s}{q}(\overline{C}|s). \label{basic1}%
\end{equation}

\end{lemma}

\begin{proof}
Let $s_{0}\in C$ and $t_{0}\in\overline{C}$ maximize the quantity $\mu
_{s}q(t\mid s)$ amongst $s\in C$ and $t\in\overline{C}$. Then $\mu_{s_{0}%
}q(t_{0}\mid s_{0})\geq\sum_{s\in C}\mu_{s}{q}(\overline{C}|s)/|S|^{2}%
\geq\zeta_{q}/|S|^{2}>\varepsilon\zeta_{q}$. Since $\widehat{q}$ is
$(\varepsilon,\beta)$-close to $q$, $\widehat{q}(t_{0}\mid s_{0})\geq
(1-\beta)q(t_{0}\mid s_{0})$. In particular,
\begin{equation}
\sum_{s\in C}\mu_{s}\widehat{q}(\overline{C}\mid s)\geq\mu_{s_{0}}\widehat
{q}(t_{0}\mid s_{0})\geq(1-\beta)\mu_{s_{0}}q(t_{0}\mid s_{0})\geq
\frac{1-\beta}{|S|^{2}}\sum_{s\in C}\mu_{s}{q}(\overline{C}|s),
\label{equ basic1}%
\end{equation}
and the left hand side inequality in (\ref{basic1}) holds.

Let $s_{1}\in C$ and $t_{1}\in\overline{C}$ maximize the quantity $\mu
_{s}\widehat{q}(t\mid s)$ amongst $s\in C$ and $t\in\overline{C}$. By
(\ref{equ basic1}), $\mu_{s_{1}}\widehat{q}(t_{1}\mid s_{1})\geq\sum_{s\in
C}\mu_{s}\widehat{q}(\overline{C}\mid s)/|S|^{2}\geq(1-\beta)\zeta_{q}%
/|S|^{4}>\varepsilon$. Since $\widehat{q}$ is $(\varepsilon,\beta)$-close to
$q$, $q(t_{1}\mid s_{1})\geq\widehat{q}(t_{1}\mid s_{1})/(1+\beta)$.
Therefore
\[
\sum_{s\in C}\mu_{s}q(\overline{C}\mid s)\geq\mu_{s_{1}}q(t_{1}\mid s_{1}%
)\geq\frac{1}{1+\beta}\mu_{s_{1}}\widehat{q}(t_{1}\mid s_{1})\geq\frac
{1}{(1+\beta)|S|^{2}}\sum_{s\in C}\mu_{s}\widehat{q}(\overline{C}\mid s),
\]
and the right hand side inequality holds as well.
\end{proof}

\begin{lemma}
\label{lemma main}Let $C\subset S$ and $s\in C$ be given. For every $g\in
G^{\beta}(C)$ (resp. $g\in\widehat{G}^{\beta}(C)$) $\mu_{s}q(g(s)\mid
s)\geq\varepsilon\zeta_{q}$ (resp. $\mu_{s}\widehat{q}(g(s)\mid s)\geq
\varepsilon\zeta_{q})$.
\end{lemma}

Note that the second claim is not symmetric to the first, since in both we use
the stationary distribution of $q$.

\begin{proof}
The proof is quite similar for $g\in G^{\beta}(C)$ and $g\in\widehat{G}%
^{\beta}(C)$. We prove the lemma for the former, and mention where the proof
for the latter differs.

Let $g \in G^{\beta}(C)$ be arbitrary. The proof is by induction over the
number of states in $C$.

If $|C|=1$, then $C=\{s\}$ for some $s\in S$. Since $g$ is $\beta$-maximal,
$\mu_{s}q(g(s)\mid s)\geq\beta/|S|\mu_{s}q(\overline{C}|s)\geq\frac{\zeta_{q}%
}{\left|  S\right|  }\beta$ (for $g\in\widehat{G}^{\beta}(C)$, by Lemma
\ref{lemmabasic1}$,$ $\mu_{s}\widehat{q}(g(s)\mid s)\geq\frac{\beta}{\left|
S\right|  }\mu_{s}\widehat{q}(\overline{C}|s)\geq\beta\frac{1-\beta}{\left|
S\right|  ^{3}}\mu_{s}q(\overline{C}|s)\geq\beta\frac{1-\beta}{\left|
S\right|  ^{3}}\zeta_{q}$).

Consider now the case $|C| > 1$.

We first assume that there are at least two edges of $g$ whose endpoints do
not belong to $C$. Let $s_{1}\neq s_{2}\in C$ Let $g_{i}$ be the restriction
of $g$ to $C\setminus\{s_{i}\}$, $i=1,2$. By \textbf{P1}, $g_{i}\in G^{\beta
}(C\setminus\{s_{i}\})$. Since any edge of $g$ is an edge of $g_{1}$ or
$g_{2}$ (or both), the induction hypothesis applied to $C\setminus\{s_{i}\}$
and $g_{i}$, $i=1,2$, implies that the claim holds for $g$.

\bigskip

Assume now that there is a unique state $s_{1}\in C$ such that $g(s_{1}%
)\not \in C$. Let $g_{1}$ be the restriction of $g$ to $C\setminus\{s_{1}\}$.
By \textbf{P1}, $g_{1}\in G^{\beta}(C\setminus\{s_{1}\})$. By the induction
hypothesis applied to $C\setminus\{s_{1}\}$ and $g_{1}$, $\mu_{s}q(g(s)\mid
s)\geq\varepsilon\zeta_{q}$ for every $s\in C\setminus\{s_{1}\}$. Thus, it
remains to show that $\mu_{s_{1}}q(g(s_{1})\mid s_{1})\geq\varepsilon\zeta
_{q}$.

Let $s_{2}\in C$ maximize the quantity $\mu_{s}q(\overline{C}\mid s)$ amongst
$s\in C$ (for $g\in\widehat{G}^{\beta}(C)$, it is chosen to maximize $\mu
_{s}\widehat{q}(\overline{C}\mid s)$). By the definition of $\zeta_{q}$,
$\mu_{s_{2}}q(\overline{C}\mid s_{2})\geq\zeta_{q}/|S|$ (for $g\in\widehat
{G}^{\beta}(C)$, by Lemma \ref{lemmabasic1}, $\mu_{s_{2}}\widehat{q}%
(\overline{C}\mid s_{2})\geq(1-\beta)\zeta_{q}/|S|^{3}$). Let $\widehat{g}\in
G^{1}(S\setminus C)$ (for $g\in\widehat{G}^{\beta}(C),$ one also chooses
$\widehat{g}\in G^{1}(S\setminus C)$) . By \textbf{P0} and \textbf{P2},
$\widehat{g}\cup g_{1}\in G(S\setminus\{s_{1}\})$.

Let $\overline{g}\in G^{1}(S\setminus\{s_{2}\})$. Since $\overline
{g}_{|S\setminus C}$ is a $S\setminus C$-graph, we have $p(\widehat{g})\geq
p(\overline{g}_{|S\setminus C})$. Since for every $t\in\overline{C}$,
$\overline{g}_{|C\setminus\{s_{2}\}}\cup(s_{2},t)$ is a $C$-graph, $p(g)\geq
p(\overline{g}_{|C\setminus\{s_{2}\}})q(t\mid s_{2})$. In particular,
$p(\widehat{g})p(g)\geq\beta p(\overline{g})q(t\mid s_{2})$ for every
$t\in\overline{C}$, and therefore $p(\widehat{g})p(g)\geq\frac{\beta}%
{|S|}p(\overline{g})q(\overline{C}\mid s_{2})$.

Denote $\sum=\sum_{y\in S}\sum_{g\in G(S\setminus\{y\})}p(g)$. By
(\ref{invariant}), $\mu_{s}=\frac{1}{\sum}\sum_{g\in G(S\setminus\{s\})}p(g)$.
In particular,
\begin{align*}
\frac{\zeta_{q}}{|S|}  &  \leq\mu_{s_{2}}q(\overline{C}\mid s_{2})\leq
\frac{\sum_{g\in G(S\setminus\{s_{2}\})}p(g)}{\sum}\times q(\overline{C}\mid
s_{2})\\
&  \leq\frac{Lp(\overline{g})q(\overline{C}\mid s_{2})}{\sum}\leq\frac
{L\times|S|}{\sum\beta}p(\widehat{g})p(g)\\
&  =\frac{L\times|S|}{\sum\beta}p(\widehat{g}\cup g\setminus(s_{1}%
,g(s_{1})))q(g(s_{1})\mid s_{1})\\
&  \leq\frac{L\times|S|}{\sum\beta}\sum_{g\in G(S\setminus\{s_{1}%
\})}p(g)\times q(g(s_{1})\mid s_{1})\\
&  =\frac{L\times|S|}{\beta}\mu_{s_{1}}q(g(s_{1})\mid s_{1}).
\end{align*}
But then $\mu_{s_{1}}q(g(s_{1})\mid s_{1})\geq\varepsilon\zeta_{q}$, as
desired. The calculation for $g\in\widehat{G}^{\beta}(C)$ is analogous.
\end{proof}

\begin{corollary}
\label{propcore} For every proper subset $C$ of $S$,%
\begin{equation}
\left|  1-\frac{\widehat{p}(g)}{p(g)}\right|  \leq(|S|+1)\beta\text{, for
every }g\in G^{\beta}(C)\cup\widehat{G}^{\beta}(C) \label{equ propcore}%
\end{equation}
and%
\[
\left|  \frac{\sum_{g\in H}\widehat{p}(g)}{\sum_{g\in H}p(g)}-1\right|
<(|S|+1)\beta\text{, where }H=G^{\beta}(C)\cup\widehat{G}^{\beta}(C).
\]

\end{corollary}

Thus, the weights of $\beta$-maximal graphs under $q$ and $\widehat{q}$ are close.

\begin{proof}
Note first that the second inequality follows immediately from the first one.
Let us prove (\ref{equ propcore}). Let $g\in G^{\beta}(C)$. By Lemma
\ref{lemma main}, $\mu_{s}q(g(s)\mid s)\geq\varepsilon\zeta_{q}$ for every
$s\in C$. Since $\widehat{q}$ is $(\varepsilon,\beta)$-close to $q$,
$(1-\beta)q(g(s)\mid s)\leq\widehat{q}(g(s)\mid s)\leq(1+\beta)q(g(s)\mid s)$.
Multiplying this inequality over $s\in C$ yields $(1-\beta)^{|C|}%
p(g)\leq\widehat{p}(g)\leq(1+\beta)^{|C|}p(g)$, and (\ref{equ propcore}) follows.

The proof for $g\in\widehat{G}^{\beta}(C)$ is similar.
\end{proof}

\subsection{Proof of Theorem \ref{theoremmain}}

\begin{proposition}
\label{prop qhatirr} The transition matrix $\widehat{q}$ is irreducible.
\end{proposition}

\begin{proof}
It is enough to prove that for every non-empty subset $C\subset S,$ there
exists $s\in C$, and $t\not \in C$ such that $\widehat{q}(t\mid s)>0$.

Let $s_{1}\in C$ and $t_{1}\not \in C$ be such that $\mu_{s_{1}}q(t_{1}\mid
s_{1})\geq\zeta_{q}/|S|^{2}>\varepsilon\zeta_{q}$. Since $\widehat{q}$ is
$(\varepsilon,\beta)$-close to $q$, $\widehat{q}(t_{1}\mid s_{1}%
)>(1-\beta)q(t\mid s)>0$.
\end{proof}

We need the following technical Lemma.

\begin{lemma}
\label{lemma technical}

\begin{enumerate}
\item Let $(a_{i})_{i=1}^{I}$ and $(b_{i})_{i=1}^{I}$ be positive numbers, and
let $\varepsilon> 0$. If $\left|  \frac{a_{i}}{b_{i}} - 1 \right|  <
\varepsilon$ for every $i=1,\ldots,I$ then $\left|  \frac{\sum_{i=1}^{I}
a_{i}}{\sum_{i=1}^{I} b_{i}} - 1 \right|  < \varepsilon$ and $\left|
\frac{\min\{a_{1},a_{2},\ldots,a_{I}\}}{\min\{b_{1},b_{2},\ldots,b_{I}\}} - 1
\right|  < \varepsilon$. \label{t1}

\item Let $\varepsilon\in(0,1/3)$, and let $a,A,b,B > 0$. If $\left|  \frac
{a}{b} - 1 \right|  < \varepsilon$ and $\left|  \frac{b}{B} - 1 \right|  <
\varepsilon$ then $\left|  \frac{a/b}{A/B} - 1 \right|  < 3\varepsilon$.
\label{t2}
\end{enumerate}
\end{lemma}

\begin{proof}
The proof of the first part is left to the reader. For the second part, note
that $1/(1+\varepsilon)<B/b<1/(1-\varepsilon)$, which implies that
$B/b-1<\varepsilon/(1-\varepsilon)$. In particular,
\[
\left|  \frac{a/b}{A/B}-1\right|  \leq\left(  \left|  \frac{a}{A}-1\right|
+1\right)  \left|  \frac{B}{b}-1\right|  +\left|  \frac{a}{A}-1\right|
<(1+\varepsilon)\frac{\varepsilon}{1-\varepsilon}+\varepsilon<3\varepsilon.
\]

\end{proof}

\begin{proposition}
\label{prop 21} For each $s\in S$,
\[
|1-\frac{\widehat{\mu}_{s}}{\mu_{s}}|<18\beta L.
\]

\end{proposition}

\begin{proof}
Fix $s\in S$. By (\ref{invariant}),
\[
\mu_{s}=\frac{\sum_{G(S\backslash\left\{  s\right\}  )}p(g)}{\sum_{y\in S}%
\sum_{G(S\backslash\left\{  y\right\}  )}p(g)}\text{ and }\widehat{\mu}%
_{s}=\frac{\sum_{G(S\backslash\left\{  s\right\}  )}\widehat{p}(g)}{\sum_{y\in
S}\sum_{G(S\backslash\left\{  y\right\}  )}\widehat{p}(g)}.
\]
For every $y\in S$, define $H_{y}=G^{\beta}(S\backslash\left\{  y\right\}
)\cup\widehat{G}^{\beta}(S\backslash\left\{  y\right\}  )$. Define%
\[
\mu_{s}^{\prime}=\frac{\sum_{H_{s}}p(g)}{\sum_{y\in S}\sum_{H_{y}}p(g)}\text{
and }\widehat{\mu}_{s}^{\prime}=\frac{\sum_{H_{s}}\widehat{p}(g)}{\sum_{y\in
S}\sum_{H_{y}}\widehat{p}(g)}.
\]
By Lemma \ref{lemma h} and Lemma \ref{lemma technical}, $\left|  \frac{\mu
_{s}}{\mu_{s}^{\prime}}-1\right|  <3\beta L$ and $\left|  \frac{\widehat{\mu
}_{s}}{\widehat{\mu}_{s}^{\prime}}-1\right|  <3\beta L$. By Lemmas
\ref{lemma h} and \ref{lemma technical}, $\left|  \frac{\widehat{\mu}%
_{s}^{\prime}}{\mu_{s}^{\prime}}-1\right|  <3(\left|  S\right|  +1)\beta$.
Since $L\geq\left|  S\right|  \geq2$, the result follows by Lemma
\ref{lemma technical}(\ref{t2}).
\end{proof}

\subsection{Proof of Theorem \ref{theorem main2}}

\begin{proposition}
\label{prop 22} For every proper subset $C$ of $S$, every $s\in C$ and
$t\not \in C$,
\[
\left|  {\mathbf{Q}}_{s,q}(t\mid C)-{\mathbf{Q}}_{s,\widehat{q}}(t\mid
C)\right|  <12\beta L.
\]

\end{proposition}

\begin{proof}
Denote $H = G^{\beta}(C) \cup\widehat G^{\beta}(C)$, and $H_{s,t} = H \cap
G_{s,t}(C)$.

Assume first that $H_{s,t}\neq\emptyset$. By (\ref{distribution}), one has
\[
\left|  \mathbf{Q}_{s,q}(t|C)-\frac{\sum_{H\cap G_{s,t}(C)}p(g)}{\sum
_{G(C)}p(g)}\right|  \leq\beta L\text{.}%
\]
Since $\left|  \frac{\sum_{H}p(g)}{\sum_{G(C)}p(g)}-1\right|  \leq\beta L$,
this yields, by Lemma \ref{lemma technical}(\ref{t2}),
\begin{equation}
\left|  \mathbf{Q}_{s,q}(t|C)-\frac{\sum_{H\cap G_{s,t}(C)}p(g)}{\sum_{H}%
p(g)}\right|  \leq\beta L+3\beta L\leq4\beta L\text{,} \label{step11}%
\end{equation}
and a similar inequality holds with $q$ replaced by $\widehat{q}$.

By Corollary \ref{propcore} and Lemma \ref{lemma technical}(\ref{t1}),%
\begin{equation}
\left|  \frac{\sum_{H\cap G_{s,t}(C)}\widehat{p}(g)}{\sum_{H\cap G_{s,t}%
(C)}p(g)}-1\right|  \leq\beta(\left|  S\right|  +1)\text{ and }\left|
\frac{\sum_{H}\widehat{p}(g)}{\sum_{H}p(g)}\right|  \leq\beta(\left|
S\right|  +1)\text{.} \label{step12}%
\end{equation}
By Lemma \ref{lemma technical}(\ref{t2}), (\ref{step12}) implies%
\[
\left|  \frac{\sum_{H\cap G_{s,t}(C)}p(g)}{\sum_{H}p(g)}-\frac{\sum_{H\cap
G_{s,t}(C)}\widehat{p}(g)}{\sum_{H}\widehat{p}(g)}\right|  \leq3\beta(\left|
S\right|  +1),
\]
which implies, using (\ref{step11}),
\[
\left|  \mathbf{Q}_{s,q}(t|C)-\mathbf{Q}_{s,\widehat{q}}(t|C)\right|
\leq\beta(8L+3(\left|  S\right|  +1)).
\]

If, on the other hand, $H_{s,t} = \emptyset$, then by (\ref{distribution}) and
the definition of $H_{s,t}$, ${\mathbf{Q}}_{s,q}(t \mid C),{\mathbf{Q}%
}_{s,\widehat q}(t \mid C) \leq\beta L $.
\end{proof}

\begin{proposition}
\label{prop 23} For every proper subset $C$ of $S$,
\[
\frac{1}{2\left|  S\right|  ^{2}}K_{C}\leq\widehat{K}_{C}\leq2\left|
S\right|  ^{2}K_{C}.
\]

\end{proposition}

\begin{proof}
We first argue that
\begin{equation}
K_{C}=\frac{\sum_{s\in C}\mu_{s}}{\sum_{s\in C}\mu_{s}q(\overline{C}\mid s)}.
\label{equ q}%
\end{equation}
Indeed, define the r.v. $\rho_{n}$ as the average length of visits to $C$ that
end before stage $n+1$:%
\begin{equation}
\rho_{n}=\frac{\sum_{p=1}^{n}1_{s_{p}\in C}}{\sum_{p=1}^{n}1_{s_{p}\in
C}1_{s_{p+1}\notin C}}. \label{equ rho}%
\end{equation}
By the ergodic theorem, the sequence $(\rho_{n})$ converges, $\mathbf{P}%
_{s_{1},q}$-a.s. to $K_{C}$, while the right hand side in (\ref{equ rho})
converges $\mathbf{P}_{s_{1},q}$-a.s. to $\frac{\sum_{s\in C}\mu_{s}}%
{\sum_{s\in C}\mu_{s}q(\overline{C}|s)}$. The identity (\ref{equ q}) follows.

By Proposition \ref{prop 21}, for every $s\in C$,
\begin{equation}
(1-18\beta L)\mu_{s}<\widehat{\mu}_{s}<(1+18\beta L)\mu_{s}. \label{equ 23.1}%
\end{equation}
By Lemma \ref{lemmabasic1},
\begin{equation}
\frac{1-\beta}{|S|^{2}}\sum_{s\in C}\mu_{s}{q}(\overline{C}|s)\leq\sum_{s\in
C}\mu_{s}\widehat{q}(\overline{C}|s)\leq(1+\beta)|S|^{2}\sum_{s\in C}\mu
_{s}{q}(\overline{C}|s). \label{equ 23.2}%
\end{equation}
Eqs. (\ref{equ 23.1}) and (\ref{equ 23.2}) yield
\begin{equation}
\frac{(1-18\beta L)(1-\beta)}{|S|^{2}}\sum_{s\in C}\mu_{s}{q}(\overline
{C}|s)\leq\sum_{s\in C}\widehat{\mu}_{s}\widehat{q}(\overline{C}%
|s)\leq(1+\beta)(1+18\beta L)|S|^{2}\sum_{s\in C}\mu_{s}{q}(\overline{C}|s).
\label{equ 23.3}%
\end{equation}
Summing up Eq. (\ref{equ 23.1}) over $s\in C$ gives
\begin{equation}
(1-18\beta L)\sum_{s\in C}\mu_{s}<\sum_{s\in C}\widehat{\mu}_{s}<(1+18\beta
L)\sum_{s\in C}\mu_{s}. \label{equ 23.4}%
\end{equation}
The Proposition follows by dividing (\ref{equ 23.4}) by (\ref{equ 23.3}).
\end{proof}

\section{Proof of the variations}

We here prove Theorems \ref{theoremvariation} and \ref{theorem variation2}. We
shall follow the previous proofs, and will point out which changes are needed.
We let $a,\varepsilon,\beta$ be given, that satisfy the assumptions of Theorem
\ref{theoremvariation}. The result of Section \ref{secestimates} still hold
for every proper subset $C$ of $S_{1}$, namely Lemmas \ref{lemmabasic1},
\ref{lemma main} and Corollary \ref{propcore} are still valid, provided the
assumption $C\subset S$ is replaced by the assumption $C\subset S_{1}$.

\subsection{Proof of Theorem \ref{theoremvariation}}

We need the following observation.

\begin{lemma}
\label{lemma a} For every $y\in S_{1}$, there exists a $(a/L)^{\left|
S\right|  }$-maximal graph $\overline{g}\in G(S_{1}\backslash y)$ such that
all paths of $\overline{g}$ lead to $y$.
\end{lemma}

\begin{proof}
By \textbf{P2}, for every $s\in S_{1}\backslash y$ there is a $\frac{a}{L}%
$-maximal $S_{1}\backslash y$-graph $g_{s}$ in which $s$ leads to a state in
$y$. Let $h_{s}$ be the path in $g_{s}$ that connects $s$ to $y\;$(this is a
set of edges). Let $\overline{g}$ be a $S_{1}\backslash y$-graph that is
contained in $\cup_{s\in S_{1}\backslash y}h_{s}$. Then $\overline{g}$
satisfies the conditions.
\end{proof}

\bigskip

We next prove the two assertions of Theorem \ref{theoremvariation}.

\begin{lemma}
All states of $S_{1}$ belong to the same recurrent set for $\widehat{q}$.
\end{lemma}

\begin{proof}
It is enough to prove that for each $C\subset S_{1}$, there exists $s\in C$
and $t\in\overline{C}$ such that $\widehat{q}(t|s)>0$. The proof of
Proposition \ref{prop qhatirr} still applies, provided $\zeta_{q}$ is replaced
by $\zeta_{q}^{1}$.
\end{proof}

\bigskip

\begin{lemma}
For each $s\in S_{1},$%
\[
\left|  1-\frac{\widehat{\mu}(s|S_{1})}{\mu(s|S_{1})}\right|  \leq18\beta L.
\]

\end{lemma}

\begin{proof}
The proof goes essentially as in Proposition \ref{prop 21}. Set $\eta
=\beta/(a/L)^{|S|}<(a/L)^{S}$, and fix $s\in S_{1}$. By (\ref{invariant}),
\[
\mu(s|S_{1})=\frac{\sum_{G(S\backslash\left\{  s\right\}  )}p(g)}{\sum_{y\in
S_{1}}\sum_{G(S\backslash\left\{  y\right\}  )}p(g)}\text{ and }\widehat{\mu
}(s|S_{1})=\frac{\sum_{G(S\backslash\left\{  s\right\}  )}\widehat{p}(g)}%
{\sum_{y\in S_{1}}\sum_{G(S\backslash\left\{  y\right\}  )}\widehat{p}(g)}.
\]
For every $y\in S_{1}$, define $H_{y}=G^{\eta}(S\backslash\left\{  y\right\}
)\cup\widehat{G}^{\eta}(S\backslash\left\{  y\right\}  )$. Define%
\[
\mu^{\prime}(s|S_{1})=\frac{\sum_{H_{s}}p(g)}{\sum_{y\in S_{1}}\sum_{H_{y}%
}p(g)}\text{ and }\widehat{\mu}^{\prime}(s|S_{1})=\frac{\sum_{H_{s}}%
\widehat{p}(g)}{\sum_{y\in S_{1}}\sum_{H_{y}}\widehat{p}(g)}.
\]
Fix for a moment $y\in S_{1}$. By Lemma \ref{lemma a} there is a $(a/L)^{|S|}%
$-maximal $S_{1}\setminus\{y\}$-graph $\overline{g}$ such that all its paths
lead to $y$. Let $g\in G^{\eta}(S\backslash\left\{  y\right\}  )$, and
$g_{S_{1}\backslash\left\{  y\right\}  }$, $g_{S\backslash S_{1}}$ its
restrictions to $S_{1}\backslash\left\{  y\right\}  $ and $S\backslash S_{1}$.
Using the above remark, the graph $\overline{g}\cup g_{S\backslash S_{1}}$ is
a $S\backslash\left\{  y\right\}  $-graph. Therefore, $g_{S_{1}\backslash
\left\{  y\right\}  }$ is $\eta(a/L)^{\left|  S\right|  }$-maximal (= $\beta
$-maximal). By Corollary \ref{propcore} one has $\left|  1-\frac{\widehat
{p}(g_{S_{1}\backslash\left\{  y\right\}  })}{p(g_{S_{1}\backslash\left\{
y\right\}  })}\right|  <(\left|  S\right|  +1)\beta$. Since $q$ and
$\widehat{q}$ coincide outside $S_{1}$, $p(g_{S\backslash S_{1}})=\widehat
{p}(g_{S\backslash S_{1}})$. Thus, $\left|  1-\frac{\widehat{p}(g)}%
{p(g)}\right|  <(\left|  S\right|  +1)\beta$. Lemma \ref{lemma h} and Lemma
\ref{lemma technical} implies that $\left|  \frac{\mu(s|S_{1})}{\mu^{\prime
}(s|S_{1})}-1\right|  <3\beta L$ and $\left|  \frac{\widehat{\mu}(s|S_{1}%
)}{\widehat{\mu}^{\prime}(s|S_{1})}-1\right|  <3\beta L$. By Lemmas
\ref{lemma h} and \ref{lemma technical}, $\left|  \frac{\widehat{\mu}^{\prime
}(s|S_{1})}{\mu^{\prime}(s|S_{1})}-1\right|  <3(\left|  S\right|  +1)\beta$.
Since $L\geq\left|  S\right|  \geq2$, the Lemma follows by Lemma
\ref{lemma technical}(\ref{t2}).
\end{proof}

\bigskip

\subsection{Proof of Theorem \ref{theorem variation2}}

\bigskip

The proof of the first two assertions in Theorem \ref{theorem variation2} is
identical to the proof of the two assertions in Theorem \ref{theorem main2}
(see Propositions \ref{prop 22} and \ref{prop 23}). We omit it. We now
prove\ a slightly strengthened version of the last assertion.

\bigskip

\begin{proposition}
Let $\eta\leq\varepsilon\zeta_{q}^{1}$ be such that$\left|  1-\frac
{\widehat{q}(t|s)}{q(t|s)}\right|  \leq\beta$ whenever $\mu_{s}\max
(q(t|s),\widehat{q}(t|s))\geq\eta$. One has
\[
\frac{1}{c}K_{S_{1}}\leq\widehat{K}_{S_{1}}\leq cK_{S_{1}}\text{ or }K_{S_{1}%
},\widehat{K}_{S_{1}}\geq\frac{1}{2\left|  S\right|  }\times\frac{\mu_{S_{1}}%
}{\eta}.
\]

\end{proposition}

The last statement in Theorem \ref{theorem variation2} corresponds to the case
$\eta=\varepsilon\zeta_{q}^{1}$.

\begin{proof}
Fix $s\in S_{1}$. By (\ref{equ q}),
\[
K_{S_{1}}=\frac{1}{\sum_{t\in S_{1}}\mu(t|S_{1})q(\overline{S_{1}}\mid t)},
\]
and a similar equality holds for $\widehat{K}_{S_{1}}$, involving
$\widehat{\mu}$ and $\widehat{q}$. By Theorem \ref{theoremvariation}%
(\ref{variantmu}) and Lemma \ref{lemma technical}, the ratio between
$\widehat{K}_{S_{1}}$ and $\frac{1}{\sum_{t\in S_{1}}\mu(t|S_{1})\widehat
{q}(\overline{S_{1}}\mid t)}$ is between $1-54\beta L$ and $1+54\beta L$.

If for every $t\in S_{1}$ and $u\not \in S_{1}$, $\mu_{t}q(u\mid t)<\eta$ and
$\mu_{t}\widehat{q}(u\mid t)<\eta$, then $K_{S_{1}}\geq\frac{\mu_{S_{1}}%
}{|S|^{2}\eta}$ and $\widehat{K}_{S_{1}}\geq(1-54\beta L)\times\frac
{\mu_{S_{1}}}{|S|^{2}\eta}$, as desired.

If, on the other hand, there exist $t\in S_{1}$ and $u\not \in S_{1}$ such
that $\mu_{t}q(u\mid t)\geq\eta$ or $\mu_{t}\widehat{q}(u\mid t)\geq\eta$ then
$|1-\frac{\widehat{q}(u\mid t)}{q(u\mid t)}|\leq\beta$, and therefore $\mu
_{t}q(u\mid t)\geq(1-\beta)\eta$ and $\mu_{t}\widehat{q}(u\mid t)\geq
(1-\beta)\eta$. For every $t\in S_{1}$ and $u\not \in S_{1}$ such that
$\mu_{t}q(u\mid t)<\eta$ and $\mu_{t}\widehat{q}(u\mid t)<\eta$ we have
$\mu_{t}q(u\mid t)\leq\sum_{t\in S_{1}}\mu_{t}q(\overline{S}_{1}\mid t)$ and
$\mu_{t}\widehat{q}(u\mid t)\leq\sum_{t\in S_{1}}\mu_{t}\widehat{q}%
(\overline{S}_{1}\mid t)$. It follows that the ratio between $\sum_{t\in
S_{1}}\mu_{t}q(\overline{S}_{1}\mid t)$ and $\sum_{t\in S_{1}}\mu_{t}%
\widehat{q}(\overline{S}_{1}\mid t)$ is at most $|S|^{2}$. The result follows.
\end{proof}

\newpage

\bibliographystyle{plain}
\bibliography{eco,jeux,math,optim,perso}

\end{document}